\documentclass[11pt]{article}

\usepackage[centertags]{amsmath}
\usepackage[dvips]{graphicx}
\usepackage{amsfonts,amssymb,amsthm,enumerate,newlfont}
\usepackage{latexsym}
\usepackage{color}
 \usepackage{hyperref}

\newcommand{\dem}{\noindent{\bf Proof. \quad}}

\newtheorem {teo} {\bf Theorem\,} [section]
\newtheorem {prop} [teo] {\bf Proposition}
\newtheorem {coro} [teo]{\bf Corollary}
\newtheorem {Lemma} [teo] {\bf Lemma}
\newtheorem {defn} [teo] {\bf Definition}

\newtheorem {Remark} [teo] {\bf Remark}

  \usepackage{paralist}
  \usepackage{graphics} 
  \usepackage{epsfig} 
\usepackage{graphicx}  \usepackage{epstopdf}

\newcommand{\R}{\mathbb{R}}

\baselineskip
\title{A GRADIENT FLOW GENERATED BY A NONLOCAL MODEL OF A NEURAL FIELD  IN AN UNBOUNDED DOMAIN}

\begin{document}

\author{\normalsize SEVERINO H. DA SILVA$^{2}$\thanks{Unidade Acad\^emica de Matem\'atica (UAMat) - Universidade Federal de Campina Grande (UFCG), Avenida
  Apr\'igio Veloso, 882, Bairro Universit\'ario, Caixa Postal: 10.044, 58109-970, Campina Grande-PB, Brazil.
  Partially supported by CAPES/CNPq-Brazil.}  \hspace{0.1cm}
 ANT\^ONIO L. PEREIRA$^{1}$\thanks{Instituto de Matem\'atica e Estat\'istica (IME)-Universidade de S\~ao Paulo (USP),
 Rua do Mat\~ao, 1010, Cidade Universit\'aria, 05508-090, S\~ao Paulo-SP, Brazil.
 Partially supported by CNPq-Brazil grants 2003/11021-7, 03/10042-0.},\hspace{0.1cm}
  \\{\scriptsize E-mail: horacio@mat.ufcg.edu.br and alpereir@ime.usp.br}}
\maketitle

\section*{Abstract}
In this paper we consider the non local evolution equation
$$
\frac{\partial u(x,t)}{\partial t} + u(x,t)=
\int_{\mathbb{R}^{N}}J(x-y)f(u(y,t))\rho(y)dy+ h(x). 
$$
 We show that this equation defines a continuous flow in both the space 
$C_{b}(\mathbb{R}^{N})$  of bounded  continuous functions and the space
 $C_{\rho}(\mathbb{R}^{N})$ of continuous functions $u$ such that $u \cdot \rho$ is bounded, where 
$\rho $ is a convenient "weight function"'.     
  We show the existence of an absorbing ball for  the flow in $C_{b}(\mathbb{R}^{N})$ and the existence of a global compact attractor for the flow 
	 in $C_{\rho}(\mathbb{R}^{N})$, under additional conditions on the nonlinearity.

We then  exhibit a continuous Lyapunov function which is well defined in the whole phase space and continuous in the 
 $C_{\rho}(\mathbb{R}^{N})$  topology,  allowing  the  characterization of the attractor as the  unstable set  of the
equilibrium point set. We also illustrate our result with a
concrete example.
\\

\noindent { \em  2010 Mathematics  Subject Classification: 45J05, 37B25.}\\

\noindent{ \it{Keywords:} Nonlocal problem, neural field, weighted space,  global
attractor, Lyapunov functional}

\section{Introduction}

We consider here the non local evolution equation
\begin{equation}
\frac{\partial u(x,t)}{\partial t} + u(x,t)=
\int_{\mathbb{R}^{N}}J(x-y)f(u(y,t))\rho(y)dy+ h(x), 
\label{eq_nonloc}
\end{equation}
where
 $f$ is  a  continuous  real function, $J: \mathbb{R}^{N} \to \mathbb{R}$ is a  non negative  integrable function,   $\rho: \mathbb{R}^{N} \to \mathbb{R}$ is a symmetric non negative bounded  "`weight"'  function with $\int_{\mathbb{R}^{N}}\rho(x)d(x) < \infty$ and   $h$ is a  bounded continuous function. Additional hypotheses will be added when needed in the sequel.


We can rewrite equation (\ref{eq_nonloc}) as
$$
\frac{\partial u(x,t)}{\partial t} + u(x,t)=
J{*}_{\rho}(f\circ u )(x,t)+ h(x), \,\,\, h  \geq 0,
$$
where the $*_{\rho}$ above denotes convolution product with respect to the measure  $d\mu(y)=\rho(y)dy$, that is
\begin{equation} \label{defJrho}
J {*}_{\rho} (v)(x):=\int_{\mathbb{R}^{N}}J(x-y) v(y))  \, d\mu(y)
= \int_{\mathbb{R}^{N}}J(x-y) v(y))  \,\rho(y)  \, d(y) .
\end{equation}

Equation (\ref{eq_nonloc}) is a variation of the equation derived by Wilson and Cowan, \cite{Wilson},
to model neuronal activity. There are also other variations of this model in the literature (see, for example, \cite{Amari}, \cite{Chen}, \cite{Silva3}, \cite{Michel} and \cite{French}).

 The function $u(x,t)$ denotes the mean membrane potential of a patch of tissue located at position $x\in
\mathbb{R}^{N}$ at time $t\geq 0$. The connection function $J(x)$ determines the coupling between the elements at position $x$ with
the element at position $y$. The  (usually non negative nondecreasing function)  $f(u)$ gives the neural firing rate, or averages rate at which
spikes are generated, corresponding to an activity level $u$. The function $h$ denotes an external stimulus applied
to the entire neural field. Let us denote  by $ S(x,t)=f(u(x,t)) $ the firing rate of a neuron at position $x$ at time $t$. 
The neurons at a point $x$ are said to be active if $S(x,t)>0$, (see \cite{Amari}, \cite{Amari2} and \cite{Rubin}).

There is already a vast literature on the analysis of similar neural field models, (see \cite{Amari}, \cite{Amari2}, \cite{Chen}, \cite{Silva},
\cite{Silva2}, \cite{Silva3}, \cite{Silva4}, \cite{Silva5}, \cite{Ermentrout},
\cite{Ermentrout2}, \cite{French}, \cite{Kishimoto} and \cite{Krisner},
\cite{Kubota}, \cite{Laing}, \cite{Rubin}).
However, their  asymptotic behavior have not
been fully analyzed in the case of unbounded domains. In particular, the "Lyapunov functional" appearing in the literature is  not well defined in the whole phase space, (see, for example, \cite{Michel} \cite{French} and \cite{Kubota}). One advantage of our model 
is that we will be able to define a continuous Lyapunov functional which is  well defined in the whole phase space, (see (\ref{Lyap_func})  in Section \ref{SecLyapunov}).

This paper is organized as follows. In Section \ref{flow_Cb}, we consider the flow generated by (\ref{eq_nonloc})  in the phase space of continuous bounded  functions. In Subsection  \ref{Well-posedness_rho}, we   prove that  the Cauchy problem for (\ref{eq_nonloc}) is well posed  in this phase space  with globally defined solutions,  and, in Subsection \ref{Absorb_sec}, we prove the existence of an absorbing set for the flow generated by (\ref{eq_nonloc}).
In Section \ref{flow_Crho}, we consider the problem  (\ref{eq_nonloc}) in the
phase space $ C_{\rho}(\mathbb{R}^{N})\equiv \{u:\mathbb{R}^{N}\rightarrow \mathbb{R}
\,\, \mbox{continuous with }\,\, \|u\|_{\rho}:= \sup_{x\in \mathbb{R}^{N}}\{|u(x)|\rho(x)\}<\infty\}
$, where $\rho$ is a convenient "weight function". In this section, to obtain well-posedness, we impose more stringent conditions on the nonlinearity than in the  previous section, (see Subsection \ref{Well-posedness}). On the other hand, we  obtain stronger results, including existence of a compact global attractor for the corresponding flow. Our proof uses adaptations of the technique  used in \cite{Silva2}, replacing the compact embedding $H^{1}([-l,l])
\hookrightarrow L^{2}([-l,l])$ by the compact embedding $C^{1}(\mathbb{R}^{N})\hookrightarrow C_{\rho}(\mathbb{R}^{N})$, (see also  \cite{Silva}, \cite{Michel},  and  \cite{Pereira} for related work).
 In Section \ref{SecLyapunov}, motivated by the energy functional
from \cite{Amari2},  \cite{Silva4}, \cite{Michel}, \cite{French}, \cite{Kubota}, and \cite{Wu}, we exhibit a continuous
Lyapunov functional for the flow generated by (\ref{eq_nonloc}),  well defined in the whole phase space $C_\rho(\mathbb{R}^{N})$, and use it to prove that the flow is gradient in the sense of \cite{Hale}. 
Finally, in Section \ref{example}, we present a concrete example to illustrate  our results.

 \section{The flow in the space $ C_{b}(\mathbb{R}^{N})$ } \label{flow_Cb}

In this section, we consider the problem  (\ref{eq_nonloc}) in the
phase space 
 $$
 C_{b}(\mathbb{R}^{N})\equiv \{u:\mathbb{R}^{N}\rightarrow \mathbb{R}
 \,\, \mbox{continuous with }\,\,    \|u\|_{\infty}:=\sup_{x\in \mathbb{R}^{N}}\{|u(x)| \}  <\infty\}. 
$$

After establishing well-posedness, we prove that a ball of appropriate radius is an absorbing set for the corresponding flow.

\subsection{Well-posedness } \label{Well-posedness_rho}

 The following estimate will be useful in the sequel. The proof is straightforward and left to the reader.
\begin{Lemma} \label{Jest}
  If  $u\in C_{b}(\mathbb{R}^{N})$ then
$\|J{*}_{\rho}u\|_{\infty}  \leq    \|J\|_{L^1(R^N)} \|\rho\|_{\infty} \|u\|_{\infty},$
where $J{*}_{\rho}u$ is given by \eqref{defJrho}.
\end{Lemma}

\vspace{3mm}
 \begin{defn} \label{loclip}
If $E$ and $F$ are   normed spaces,
we say that a function  $F : E \to F$ is
   \emph{locally Lipschitz
 continuous  (or simply locally Lipschitz) } if,\footnote{}
 for any $ x_0 \in   E$, there exists  a constant
 $C$ and a ball $B = \{ x \in E \ : \
  \|x-x_0  \|< b \}$ such that,
 if $ x$ and $ y$ belong to $B$   then
 $ \|F(x) - F(y) \| \leq  C \| x-y \| $; we say that $F$ is
   \emph{ Lipschitz
 continuous on bounded sets } if the ball
 $B$ in the previous definition can chosen as any bounded ball in
$E$.
 \end{defn}

\begin{Remark}
  The two definitions in (\ref{loclip})  are equivalent if the normed space $ E  $ is locally compact.
\end{Remark}

\begin{prop}
If   $f $  is continuous,  then the map
$F:C_{b}(\mathbb{R}^{N})\rightarrow C_{b}(\mathbb{R}^{N})$, given by 
$$
F(u)=-u+J{*}_{\rho}(f\circ u)+h, 
$$
is  well defined. If $f$  is locally Lipschitz,   then $F$  Lipschitz in bounded sets.\label{Pro 2.2}
\end{prop}
\dem The first assertion is immediate. Now, from triangle inequality and Lemma \ref{Jest}, it follows that
\begin{eqnarray*}
\|F(u)-F(v)\|_{\infty} &\leq & \|v-u\|_{\infty}+\|J{*}_{\rho}(f\circ
u)-J{*}_{\rho}(f\circ
v)\|_{\infty} \\
&\leq& \|v-u\|_{\infty} +   \|J\|_{L^1(R^N)} \|\rho\|_{\infty} \|(f\circ u)-(f\circ
v)\|_{\infty}  .
\end{eqnarray*}
If $ \|u\|_{\infty},  \|v \|_{\infty} \leq R $ then $|(f\circ u)(x)-(f\circ v)(x)|
\leq k_R |u(x) - v(x)| $, where $k_R$ is a  Lipschitz constant for
$f$ in the interval $ [-R, R]$.
 It follows that

\begin{eqnarray*}
\|F(u)-F(v)\|_{\infty} \leq (1+ k_{R} \|J\|_{L^1(R^N)} \|\rho\|_{\infty})\|u-v\|_{\infty}.
\end{eqnarray*}
which concludes the proof.
\qed

\begin{teo}
 If $f$ is locally Lipschitzian, the Cauchy problem for
(\ref{eq_nonloc}) is well posed in $C_{b}(\mathbb{R}^{N})$ with globally
defined solutions. 
\end{teo}

\dem
 It follows from  Proposition \ref{Pro 2.2} and well-known results (see \cite{Daleckii} or \cite{Henry}, Theorems 3.3.3 and
3.3.4).
\qed

\subsection{Existence of an absorbing set} \label{Absorb_sec}
 
In this section, we  denote by $T(t)$ the flow generated by
(\ref{eq_nonloc}) in   $C_{b}(\mathbb{R}^{N})$. Under some additional hypotheses on the nonlinearity, we prove here the existence of an absorbing   bounded ball
 ${\cal B} \subset C_{b}(\mathbb{R}^{N})$  for $T(t)$.

We recall that a set ${\cal B} \subset C_{b}(\mathbb{R}^{N})$ is an
\emph{absorbing set}  for the flow $T(t)$ if, for any bounded set $C
\subset C_{b}(\mathbb{R}^{N})$, there is a $t_{1}=t_1(C)>0$ such that
$T(t)C \subset {\cal B}$ for any $t\geq t_{1}$, (see \cite{Teman}).

\begin{Lemma} 
Suppose that  $f$ is locally Lipschitz and  
satisfies the dissipative condition

\begin{equation}\label{dissip1}
  |f(x) | \leq \eta|x| + K, \textrm{ for any }  x\in \R.
 \end{equation}
 with $\eta \|J\|_{L^1(R^N)} \|\rho\|_{\infty} <1$. Then, if  $\eta \|J\|_{L^1(R^N)} \|\rho\|_{\infty} < \delta <1$.
the ball in  $C_b(\R^N)$, centered
at the origin  with radius
$ R = \frac{\|J\|_{L^1(R^N)} \|\rho\|_{\infty} K +  \|h \|_{\infty} }{  \delta - \|J\|_{L^1(R^N)} \|\rho\|_{\infty} \eta} $,  is  an absorbing set for the flow $T(t)$.\label{Lemma 3.1}
\end{Lemma}
\dem
Let $u(x,t)$ be the solution of
(\ref{eq_nonloc}) with initial condition $u(\cdot,0)=u_{0}$. Then, by the
variation of constants formula,
$$
u(x,t)=e^{-t}u_{0}(x)+\int_{0}^{t}e^{s-t}[J{*}_{\rho}(f\circ u)(x,s)+h]ds.
$$
From (\ref{dissip1}), there exists a constant $K$ such that
 $|f(x) | \leq \eta|x| + K $, for any $x\in \R$.

 Hence, using Lemma \ref{Jest} and \eqref{dissip1}, we obtain 

\begin{eqnarray*}
|u(x,t)| &\leq&
e^{-t} |u_{0}(x) |+\int_{0}^{t}e^{s-t}[|J{*}_{\rho}(f\circ
  u)(x,s)|+ |h(x) |]ds \\ 
 &\leq & e^{-t} \|u_0 \|_{\infty} +\int_{0}^{t}e^{s-t}[\|J{*}_{\rho}(f\circ
  u)(\cdot,s)\|_{\infty} + \|h \|_{\infty} ]ds \\
 &\leq & e^{-t} \|u_0 \|_{\infty} +\int_{0}^{t}e^{s-t}[\|J\|_{L^1(R^N)} \|\rho\|_{\infty} \|(f\circ
   u)(\cdot,s)\|_{\infty} + \|h \|_{\infty} ]ds  \\
 &\leq & e^{-t} \|u_0 \|_{\infty} +\int_{0}^{t}e^{s-t}[ \|J\|_{L^1(R^N)} \|\rho\|_{\infty} \eta \|
   u(\cdot,s)\|_{\infty} + \|J\|_{L^1(R^N)} \|\rho\|_{\infty} K +  \|h \|_{\infty} ]ds.
\end{eqnarray*}

Suppose   $  \|u(\cdot,s)\|_{\infty} \geq
\frac{1}{  \delta  - \|J\|_{L^1(R^N)} \|\rho\|_{\infty} \eta}  \left( \|J\|_{L^1(R^N)} \|\rho\|_{\infty} K +  \|h \|_{\infty} \right) $, for $0 \leq t \leq T$.
 Then,  for $t\in [0,T] $,    we obtain

  \[ e^{t} |u(x,t) |  \leq  \|u_0 \|_{\infty} 
  + \delta   \int_{0}^{t}e^{s}  \|
   u(\cdot,s)\|_{\infty} \,  ds \quad \textrm{ for any } x \in \R^N.  \]
 Taking the supremum on the left side, it follows that   
 
   \[ e^{t} \|u(x,\cdot) \|_{\infty} \, ds  \leq  \|u_0 \|_{\infty} 
  + \delta   \int_{0}^{t}e^{s}  \|
   u(\cdot,s)\|_{\infty}.
   \]
    
	 From Gronwall´s inequality, it then follows that 
	 $ e^{t}\|u( \cdot ,t)\|_{\infty} \leq    \|u_0 \|_{\infty}  e^{\delta t}  $ and, therefore  
	
	 \begin{equation} \label{exp_decay}
	  \|u( \cdot ,t)\|_{\infty} \leq    \|u_0 \|_{\infty}  e^{(\delta-1) t} ,
         \mbox{  \textrm{ for } }  t \in [0,T].
         \end{equation}

          It follows that there exists  
$ T_0 \leq  \frac{1}{(1-\delta)} \ln  \left(  \frac{\|J\|_{L^1(R^N)} \|\rho\|_{\infty} K +  \|h \|_{\infty} }{ \|u_0 \|_{\infty}( \delta  - \|J\|_{L^1(R^N)} \|\rho\|_{\infty} \eta )} \right)   $
 such that 

$$  \|u( \cdot ,T_0)\|_{\infty} \leq   
                \frac{ \|J\|_{L^1(R^N)} \|\rho\|_{\infty} K +  \|h \|_{\infty} }{ \delta- \|J\|_{L^1(R^N)} \|\rho\|_{\infty} \eta  }.$$

Also, we must have  $ \|u( \cdot ,t)\|_{\infty} \leq   
\frac{ \|J\|_{L^1(R^N)} \|\rho\|_{\infty} K +  \|h \|_{\infty} }{ \delta- \|J\|_{L^1(R^N)} \|\rho\|_{\infty} \eta  }$,
  for
  any $t\geq T_0$, since   $ \|u( \cdot ,t)\|_{\infty} $ decreases (exponentially)
  if the opposite inequality holds, by \eqref{exp_decay}.  \qed

\begin{Remark}
From \eqref{exp_decay}, it follows that the ball  $B(0,
R^{\prime})$ is positively invariant under the  flow $T(t)$ if $R^{\prime} \geq R$.
\end{Remark}

\section{The flow in the space $ C_{\rho}(\mathbb{R}^{N})$ } \label{flow_Crho}

In this section, we consider the problem  (\ref{eq_nonloc}) in the
phase space

$$
C_{\rho}(\mathbb{R}^{N})\equiv \{u:\mathbb{R}^{N}\rightarrow \mathbb{R}
\,\, \mbox{continuous with }\,\, \|u\|_{\rho}:= \sup_{x\in \mathbb{R}^{N}}\{|u(x)|\rho(x)\}<\infty\}.
$$
 
 We will need to impose more stringent conditions on the nonlinearity than in the  previous section, to obtain well-posedness. On the other hand, we will obtain stronger results, including existence of a compact global attractor for the corresponding flow.

\subsection{Well-posedness } \label{Well-posedness}

 The following result is the analogous of Lemma \ref{Jest}. The proof is again  straightforward and left to the reader.
 
\begin{Lemma} \label{Jest_rho} 
 \item  If  $u\in C_{\rho}(\mathbb{R}^{N})$ then
$\|J{*}_{\rho}u\|_{\rho} \leq  \|J \|_{L^1(\R^N)} \| \rho \|_{\infty} \|u\|_{\rho}.$
 \end{Lemma}

 \begin{prop}
If   $f $  is globally Lipschitzian,  then the map
$F:C_{b}(\mathbb{R}^{N})\rightarrow C_{b}(\mathbb{R}^{N})$, given by 
$$
F(u)=-u+J{*}_{\rho}(f\circ u)+h, 
$$
is  well defined and globally Lipschitzian.
\end{prop}
\dem  Suppose $ |f(x) - f(y)| \leq k  |x-y|$, for any $x,y \in \R$.  Then,
 in particular, $ |f(x) | \leq k  |x| + M $, where $M= f(0)$  for any $x \in \R$.
 It follows that  
  $ \|f\circ u \|_{\rho} \leq k \| u \|_{\rho} + M \|\rho\|_{\infty} $.  From Lemma
   \ref{Jest_rho}, we then obtain  
   \begin{eqnarray*} \|F(u) \|_{\rho} &\leq & \|u\|_{\rho}+\|J{*}_{\rho}(f\circ
u)\|_{\rho} \\
&\leq& \|u\|_{\rho} +   \|J \|_{L^1(\R^N)} \| \rho \|_{\infty} \|f\circ u\|_{\rho} \\
&\leq& \|u\|_{\rho} +   \|J \|_{L^1(\R^N)} \| \rho \|_{\infty} 
(  k \| u \|_{\rho} + M \|\rho\|_{\infty} ), \\
\end{eqnarray*}
so $F$ is well defined. 
   Furthermore
\begin{eqnarray*}
\|F(u)-F(v)\|_{\rho} &\leq & \|u-v\|_{\rho}+\|J{*}_{\rho}(f\circ
u)-J{*}_{\rho}(f\circ
v)\|_{\rho} \\
&\leq& \|u-v\|_{\rho} +  \|J\|_{L^1(R^N)}  \|\rho\|_{\infty}\|(f\circ u)-(f\circ
v)\|_{\rho}  \\
&\leq& \|u-v\|_{\rho } +  \|J\|_{L^1(R^N)}  \|\rho\|_{\infty}  k \| u - v\|_{\rho}  \\
 &= & ( 1+ k \|J\|_{L^1(R^N)}  \|\rho\|_{\infty})\| u - v\|_{\rho} 
  \end{eqnarray*}
Therefore $F$ is globally  Lipschitz in $C_{\rho}(\mathbb{R}^{N})$.
\qed

\begin{teo}
 If $f$ is globally  Lipschitzian, the Cauchy problem for
(\ref{eq_nonloc}) is well posed in $C_{\rho}(\mathbb{R}^{N})$ with globally
defined solutions.  
\end{teo}

\dem
 It follows from  Proposition \ref{Pro 2.2} and well-known results (see \cite{Daleckii} or \cite{Henry}, Theorems 3.3.3 and
3.3.4).
\qed

\vspace{3mm}

\subsection{Existence of an absorbing set}

In this section, we  denote by $T(t)$ the flow generated by
(\ref{eq_nonloc}) in   $C_{\rho}(\mathbb{R}^{N})$. Under some additional hypotheses on the nonlinearity, we  prove the existence of a  bounded 
ball  ${\cal B} \subset C_{\rho}(\mathbb{R}^{N})$ which is an
absorbing set for $T(t)$.

\begin{Lemma} 
Suppose that  $f$ is globally Lipschitz and  
satisfies the dissipative condition

\begin{equation}\label{dissip1_rho}
  |f(x) | \leq \eta|x| + K, \textrm{ for any }  x\in \R.
 \end{equation}
 with $ \|J\|_{L^1(\R^N)} \|\rho\|_{\infty}  \eta   <1$. Then, if  $\|J\|_{L^1(\R^N)} \|\rho\|_{\infty}  \eta  < \delta <1$,
the ball in  $C_{\rho}(\R^N)$, centered
at the origin  with radius
$ R =  \frac{\|J\|_{L^1(\R^N)} \|\rho\|_{\infty}  K +  \|h \|_{\rho}  }{ \delta- \|J\|_{L^1(\R^N)} \|\rho\|_{\infty}  \eta  }$,  is  an absorbing set for the flow $T(t)$.
\end{Lemma}
\dem
Let $u(x,t)$ be the solution of
(\ref{eq_nonloc}) with initial condition $u(\cdot,0)=u_{0}$. Then, by the
variation of constants formula,
$$
u(x,t)=e^{-t}u_{0}(x)+\int_{0}^{t}e^{s-t}[J{*}_{\rho}(f\circ u)(x,s)+h]ds.
$$
From \eqref{dissip1_rho} and  Lemma  \ref{Jest_rho},  we obtain 

\begin{eqnarray*}
|u(x,t) \rho(x)  | &\leq&
e^{-t} |u_{0}(x) \rho(x) |+\int_{0}^{t}e^{s-t}[|J{*}_{\rho}(f\circ
  u)(x,s) \rho(x)|+ |h(x) \rho(x) |]ds \\ 
  &\leq & e^{-t} \|u_0 \|_{\rho} +\int_{0}^{t}e^{s-t}[\|J{*}_{\rho}(f\circ
  u)(\cdot,s)\|_{\rho} + \|h \|_{\rho} ]ds \\
 &\leq & e^{-t} \|u_0 \|_{\rho} +\int_{0}^{t}e^{s-t}[\|J\|_{L^1(\R^N)} \|\rho\|_{\infty}
  \|(f\circ
   u)(\cdot,s)\|_{\rho} + \|h \|_{\rho} ]ds  \\
 &\leq & e^{-t} \|u_0 \|_{\rho } +\int_{0}^{t}e^{s-t}[ \|J\|_{L^1(\R^N)} \|\rho\|_{\infty}
  \eta \|
   u(\cdot,s)\|_{\rho} +  \|J\|_{L^1(\R^N)} \|\rho\|_{\infty}  K +  \|h \|_{\rho} ]ds.
\end{eqnarray*}

Suppose   $  \|u(\cdot,s)\|_{\infty} \geq
\frac{ \|J\|_{L^1(\R^N)} \|\rho\|_{\infty} K +  \|h \|_{\rho}}{  \delta  - \|J\|_{L^1(\R^N)} \|\rho\|_{\infty} \eta} 
 $, for $0 \leq t \leq T$.
 Then for $t\in [0,T] $,    we obtain

  \[ e^{t} |u(x,t) \rho(x) |  \leq  \|u_0 \|_{\rho} 
  + \delta   \int_{0}^{t}e^{s}  \|
   u(\cdot,s)\|_{\rho} \,  ds \quad \textrm{ for any } x \in \R^N.  \]
 Taking the supremum on the left side, it follows that   
 
   \[ e^{t} \|u(x,\cdot) \|_{\rho} \, ds  \leq  \|u_0 \|_{\rho} 
  + \delta   \int_{0}^{t}e^{s}  \|
   u(\cdot,s)\|_{\rho}.
   \]
    
	 From Gronwall´s inequality, it then follows that 
	 $ e^{t}\|u( \cdot ,t)\|_{\rho} \leq    \|u_0 \|_{\rho}  e^{\delta t}  $ and hence  
	
	 \begin{equation} \label{exp_decay_rho}
	  \|u( \cdot ,t)\|_{\infty} \leq    \|u_0 \|_{\infty}  e^{(\delta-1) t} ,
         \mbox{  \textrm{ for } }  t \in [0,T].
         \end{equation}

          Therefore, there exists  
$ T_0 \leq  \frac{1}{(1-\delta)} \ln  
\left(   \|J\|_{L^1(\R^N)} \|\rho\|_{\infty} K +  \|h \|_{\infty} { \|u_0 \|_{\infty}( \delta  - \|J\|_{L^1(R^N)} \|\rho\|_{\infty} \eta )} \right)   $
 such that 

$$  \|u( \cdot ,T_0)\|_{\rho} \leq   
                \frac{\|J\|_{L^1(\R^N)} \|\rho\|_{\infty}  K +  \|h \|_{\rho}  }{ \delta- \|J\|_{L^1(\R^N)} \|\rho\|_{\infty}  \eta  } $$.

Also, we must have  $  \|u( \cdot ,T_0)\|_{\rho} \leq   
                \frac{\|J\|_{L^1(\R^N)} \|\rho\|_{\infty}  K +  \|h \|_{\rho}  }{ \delta- \|J\|_{L^1(\R^N)} \|\rho\|_{\infty}  \eta  } $, for
  any $t\geq T_0$, since   $ \|u( \cdot ,t)\|_{\rho} $ decreases (exponentially)
  if the opposite inequality holds by \eqref{exp_decay_rho}.  \qed

\begin{Remark} \label{invariant_ball} 
From \eqref{exp_decay_rho}, it follows that the ball  $B_{\rho}(0,
R^{\prime})$ of radius $R^{\prime}$ in $C_{\rho}(\R^N)$ is positively invariant under the  flow $T(t)$ if $R^{\prime} \geq R$.
\end{Remark}

\subsection{Existence of a global attractor}

We denote below by  $ C_b^1(\R^N) $, the subspace of functions in  $ C_b(\R^N) $
with bounded derivatives.

\begin{Lemma} \label{compact_embed}
   The inclusion map   $i: C_b^1(\R^N) \to  C_{\rho}(\R^N)$ is compact.  

\end{Lemma}
  \dem   Let $C$ be a bounded set in   $C_b^1(\R^N)$. For any
   $l>0$, let  $\varphi: \R^N \to [0,1] $ be a smooth function satisfying
$$
\varphi(x)=\left\{\begin{array}{ccccc}
0,  \,\,if \,\, \|x\|\geq l,\\
1,\,\, if \,\,  \|x\| \leq \frac{l}{2}.\\
 
\end{array}\right.
$$
Let
$C^0( B_l) $  denote the space of continuous functions defined in
 the ball of   $\R^N$
with radius $l$ and center at the origin,  which vanish at the boundary. 
Consider the subset  $C_l$
of functions   in   $C^0( B_l) $  defined by
$$ C_l:= \{  \varphi u_{|B_l } \ \textrm{ with }  u \in C \}.  $$

Then $C_l$ is a bounded subset of $C_b^1(B_l)$  and, therefore, a 
precompact subset of $C^0( B_l) $,
 by the Arzelá-Ascoli theorem. 
 Let now ${E}_1$ be the subset of $C_{\rho}(\R^N)$ given by

 $$ {G}_l:= \{  E(u)   \ \textrm{ with }  u \in C_l \}.  $$
 where $E(u)$ is the extension   by zero outside  $ B_l$.
 Since $E$ is continuous as an operator from
 $C^0( B_l) $ into $C_{\rho}(\R^N)$, it follows that 
  $ \overline{C_l} $ is  a compact subset of $C_{\rho}(\R^N)$.  
   
 Let now 

 $$  {G^c}_l := \{ (1- \varphi)  u  \ \textrm{ with }  u \in C \}.  $$
     Let $R$ be such that  $ \|u\|_{\infty} \leq R$, for any $u \in C$. Then, for any $\epsilon> 0$, we may find $l$ such that $ 0< \rho(x) < \frac{\epsilon}{R}$,
     if $ \|x\| \geq l/2$. Then, it follows that
     $\|u\|_{\rho} \leq \epsilon$, for any $ u \in G_l^c$, that is,
     $ {G^c}_l$ is contained in the ball of radius $\epsilon$ around the origin.

     Since   ${G}_l$ is precompact, it can be covered by a finite number of balls of
     radius $\epsilon$.  Since any function $u$ in  $ C$  can be written as
      $u = u_1 + u_2$, with 
     $ u_1 = \varphi u \in G_l$ and $u_2 =   (1-\varphi) u \in G_l^c $, it follows that
      $C$ can be covered by a finite number of balls with radius $ 2 \epsilon$, for any $\epsilon > 0$. Thus $C$ is precompact as a subset of  $C_{\rho}(\R^N)$. \qed


\begin{Lemma} \label{compact_measure}  
In addition the hypotheses of Lemma \ref{dissip1_rho}, suppose that $f: \R \to \R$ is bounded and $h$ has bounded derivative. Let $C$ be a bounded set in  $C_{\rho}(\R^N)$  Then for any 
$\eta>0$, there exists $t_{\eta}$ such
that $T(t_{n}) C $, has a finite covering by balls  with radius smaller than
$\eta$.  
\end{Lemma}
\dem Let
$u(x,t)$ be the solution of (\ref{eq_nonloc}) with initial condition 
 $u_{0}\in C$. We may suppose that $C$ is contained in the ball $ B_R$ of radius $R$, centered at the origin.
By the variation of constants formula
$$
T(t)u_{0}(x)=e^{-t}u_{0}(x)+\int_{0}^{t}e^{s-t}[J{*}_{\rho}(f\circ
u)(x,s)+h(x)]ds.
$$
Write 
$$(T_{1}(t)u_{0})(x)=e^{-t}u_{0}(x)$$ 
and
$$(T_{2}(t)u_{0})(x)=\int_{0}^{t}e^{-(t-s)}[J{*}_{\rho}(f\circ u)(x,s)+h(x)]ds.$$

Let $\eta >0$ given. Then there exists $t(\eta)>0$, uniform for $u_0 \in C$, such that if
$t\geq t(\eta)$ then $\|T_{1}(t)u_{0}\|_{\rho} \leq \frac{\eta}{2}$.
In fact,
$$
|(T_{1}(t)u_{0})(x)|\rho(x)=e^{-t}|u_{0}(x)|\rho(x).
$$
Thus
$$
\|T_{1}(t)u_{0}\|_{\rho}=e^{-t}\|u_{0}\|_{\rho}.
$$
Hence, for $t > t_{\eta}=\ln \left( \frac{2 R }{\eta} \right)$, we have
$\|T_{1}(t)u_{0}\|_{\rho}\leq \frac{\eta}{2}$, for any $u_0 \in C$, that is,
$ T_{1}(t) C$ is contained in the ball of radius $ \frac{\eta}{2}$ around the origin.  


We now show that
$T_{2}(t) C_{\rho}(\R^N) $ lies in a bounded  ball  of $ C_b^{1}(\mathbb{R}^{N})$. 


In fact,  using  Lemma \ref{Jest}  we have, for any
 $u_0 \in C_{\rho}(\mathbb{R}^{N})$, 

\begin{eqnarray*}
\|T_{2}(t)u_{0}\|_{\infty}&\leq&\int_{0}^{t}e^{s-t}[\|J{*}_{\rho}(f\circ u)(\cdot, s)\|_{\infty} +\|h\|_{\infty}]ds\\
&\leq& \int_{0}^{t}e^{s-t}[\|J\|_{L^1(R^N)} \|\rho\|_{\infty}\|(f\circ u)(\cdot, s)\|_{\infty} +\|h\|_{\infty}]ds\\
&\leq&(M\|J\|_{L^1(R^N)} \|\rho\|_{\infty} + \|h\|_{\infty})\int_{0}^{t}e^{s-t}ds\\
&\leq& M\|J\|_{L^1(R^N)} \|\rho\|_{\infty} + \|h\|_{\infty},
\end{eqnarray*}

 where  $M=\|f \|_{\infty} < \infty$, and

 \begin{eqnarray*}
   \left\|\frac{\partial }{\partial x}T_{2}(t)u_{0}\right\|_{\infty}
&\leq&\int_{0}^{t}e^{s-t}[\|J'{*}_{\rho}(f\circ u)(\cdot, s)\|_{\infty} +\|h'\|_{\infty}]ds\\
&\leq& \int_{0}^{t}e^{s-t}[\|J' * \rho\|_{\infty}\|(f\circ u)(\cdot, s)\|_{\infty} +\|h'\|_{\infty}]ds\\
&\leq&(M\|J'\|_{L^1(\R^N)} \|\rho\|_{\infty} + \|h'\|_{\infty})\int_{0}^{t}e^{s-t}ds\\
&\leq& M    \|J'\|_{L^1(\R^N)} \|\rho\|_{\infty}      + \|h'\|_{\infty}
\end{eqnarray*}

 Then, for $t\geq 0$ and any $u_0 \in C_{\rho}(\R^N),$
 $\| \frac{\partial }{\partial
x}T_{2}(t)u_{0} \|_{\rho}$ is bounded by a constant independent of $t$ and $u$.


 Therefore, by Lemma  \ref{compact_embed}, it 
  follows that $\{T_{2}(t) \} C_{\rho}(\R^N),$   is compact as a subset of 
$C_{\rho}(\R^N)$ and, therefore it can be covered by a finite number of balls with radius $\frac{\eta}{2}$.

Therefore, since
$$
T(t)C =T_{1}(t) C + T_{2}(t) C, 
$$ we obtain that $ T(t)C $, can be covered by a finite number of balls of 
radius $\eta$, as claimed.
 \qed\\


In what follows we denote by $\omega(B_{\rho}(0,R))$ the
$\omega$-limit set of the ball $B_{\rho}(0,R)$.

Then as consequence from Lemma \ref{compact_measure} we have the following
result:
\begin{teo} \label{attract_rho}
Assume the same hypotheses of  Lemma \ref{compact_measure}. Then ${\cal
A}=\omega(B_{\rho}(0,R))$, is a global attractor for the flow
$T(t)$ generated by (\ref{eq_nonloc}) in $B_{\rho}(0,R)$ which is
contained in the ball of radius $B_{\rho}(0,R)$.\label{Teorema 3.3}
\end{teo}
\dem From Lemma  \ref{compact_measure}, it follows that, for any
$\eta > 0$, there exists $t_{\eta}> 0 $ such that $T(t_{\eta}) B_{\rho}(0,R)$ can be covered by a finite number of ball of radius $\eta$.  Since  $B_{\rho}(0,R)$ is
positively invariant, (see Remark \ref{invariant_ball}) we have, for any $t  \geq t_{\eta}$, $T(t)  B_{\rho}(0,R) = T(t_{\eta})T(t-t_{\eta}) B_{\rho}(0,R) \subset
 T(t_{\eta}) B_{\rho}(0,R)  $ and thus,
 $ \cup_{t \geq t_{\eta}}  T(t) B_{\rho}(0,R) \subset
 T(t_{\eta}) B_{\rho}(0,R)   $, can also be covered by a finite number of ball with radius $\eta$. 

 Therefore
 $$ {\cal A}:= \omega ( B_{\rho}(0,R)) = \cap_{t_0 \geq 0} \overline{ \cup_{t \geq t_0}  T(t) B_{\rho}(0,R)}= \cap_{t_0 \geq 0}   \overline{ T(t) B_{\rho}(0,R)},   $$
 can be covered by a finite number of balls of radius arbitrarily small radius and is closed,  so  it is a compact set.
 From the positive invariance of $B{\rho}(0,R)$ (Lemma \ref{Lemma 3.1}), it is clear that
 ${\cal A} \subset B_{\rho}(0,R)$. 
 
  It remains to prove that ${\cal A}$  attracts bounded sets of  $C_{\rho} (\R^N) $.
  It is enough to prove that it attracts the ball $ B_{\rho}(0,R)$.  
   Suppose, for contradiction, that there exist $\epsilon >0$ and  sequences $t_n \to \infty$, 
    $x_n \in  B_{\rho}(0,R)$, with $ d(T(t_n) (x_n),  {\cal A}) > \epsilon$. 

 Now,   the set $\{  T(t_n) (x_n) \ : \ n\geq n_0  \}$ is contained in   
 $ T(t_{n_0}) B_{\rho}(0,R) $, Thus for, any $\eta>0$,  it can be covered by balls with  radius $\eta$  if $n_0$ is big enough. Since the remainder of the sequence is a finite set, the same happens with the whole sequence.  It follows that the sequence
  $\{  T(t_n) (x_n) \ : \  n \in \mathbb{N} \}$ is a precompact set and  so, passing to a subsequence, it converges to a point $x_0  \in B_{\rho}(0,R)$. But then $x_0$ must belong to 
   ${\cal A}=\omega(B_{\rho}(0,R))$ and we reach a contradiction.
   
%
 This concludes the proof. \qed

\section{Existence of a Lyapunov functional } \label{SecLyapunov}

   Energy-like Lyapunov functional for models of neural fields are  well known in the literature, (see for example, \cite{Amari2}, \cite{Silva4}, \cite{Silva5}, \cite{Michel}, \cite{French}, \cite{Kubota} and \cite{Wu}. However, when  dealing with  unbounded domains, these functionals are frequently not well defined in the whole fase space, since they can assume the value  $\infty,$ at some points (see, for example, \cite{Michel}, \cite{Kubota}).

In this section, under appropriate assumptions on $f$,  we exhibit a continuous Lyapunov functional for the flow of (\ref{eq_nonloc}), which is well defined in the whole phase space $C_{\rho}(\mathbb{R}^{N})$, and used it to prove that this flow has the gradient property, in the sense of \cite{Hale}.

 Suppose that $f$ is strictly increasing. Motivated by
the  energy functionals appearing in  \cite{Amari2}, \cite{French},  \cite{Kubota}, and \cite{Wu},  we define the functional $F:C_{\rho}(\mathbb{R}^{N}) \rightarrow \mathbb{R}$  by
\begin{equation} \label{Lyap_func}
F(u)=\int_{\mathbb{R}^{N}}\left[-\frac{1}{2}f(u(x))\int_{\mathbb{R}^{N}}J(x-y)f(u(y)) \rho(y) dy +\int_{0}^{f(u(x))}f^{-1}(r)dr -h f(u(x))\right] \rho(x) dx.
\end{equation}
Equivalently, with $d\mu(x)=\rho(x)dx$, we can rewrite (\ref{Lyap_func}) as 
$$
F(u)=\int_{\mathbb{R}^{N}}\left[-\frac{1}{2}f(u(x))\int_{\mathbb{R}^{N}}J(x-y)f(u(y)) d\mu(y) +\int_{0}^{f(u(x))}f^{-1}(r)dr -h f(u(x))\right]  d\mu(x).
$$


We can then prove the following result:

\begin{prop}
  In addition to the hypotheses of Theorem \ref{attract_rho}, assume that  $f: \R \to \R$ is strictly increasing.  Then the  functional given in
    \eqref{Lyap_func}  satisfies $|F(u)|<\infty$, for all $u\in C_{\rho}(\mathbb{R}^{N})$.\label{lower}
\end{prop}
\dem We start by  noting that
$$
F(u)=F_{1}(u)+F_{2}(u)-F_{3}(u),
$$
where
$$
F_{1}(u)= -\frac{1}{2} \int_{\mathbb{R}^{N}} \int_{\mathbb{R}^{N}}\  f(u(x))J(x-y)f(u(y)) \rho(y) \rho(x) dy dx,
$$
$$
F_{2}(u)=\int_{\mathbb{R}^{N}}\left[\int_{0}^{f(u(x))}f^{-1}(r)dr\right]\rho(x)dx
$$
and
$$
F_{3}(u)=\int_{\mathbb{R}^{N}}hf(u(x))\rho(x)dx.
$$

 Let  \begin{equation}
  \label{IntegF1}
G_1(x,y) :=   f(u(x))J(x-y)f(u(y)) \rho(y) \rho(x) \end{equation}
 denote the integrand of $F_1(u) $ 
 Then, since $ M = \|f  \circ  u  \|_{\infty} <  \infty $, we obtain  

 \[
 |G_1(x,y) | \leq  M^2 J(x-y) \rho(y) \rho(x)  
   \]

and, therefore     

\begin{eqnarray} \label{boundF1}
|F_{1}(u)| 
&\leq &\frac{1}{2}\int_{\mathbb{R}^{N}}    \int_{\mathbb{R}^{N}}  M^2 J(x-y) \rho(y) \rho(x) dy dx \nonumber \\
&\leq & \frac{1}{2} {M^{2} \|J\|_{L^1(R^N)} \|\rho\|_{\infty} }  \int_{\mathbb{R}^{N}}\rho(x)dx \nonumber \\
& \leq &\frac{1}{2} {M^{2} \|J\|_{L^1(R^N)} \|\rho\|_{\infty} }   \|\rho\|_{L^1(R^N)}     ,
\end{eqnarray}

Let now  \begin{equation}
\label{IntegF2} 
G_2(x) :=   \int_{0}^{f(u(x))}f^{-1}(r)dr \rho(x) 
\end{equation}
denote the integrand of $F_2(u)$. 
 Then,

 \[  |G_2(x) | \leq  \int_{0}^{M } |f^{-1}(r)| \, dr \rho(x)
 \]
and
\begin{eqnarray} \label{boundF2}
|F_{2}(u)|&\leq& \int_{\mathbb{R}^{N}} \left[ \int_{0}^{M } |f^{-1}(r)| \, dr  \right]    \rho(x)dx \nonumber \\
&\leq&\int_{\mathbb{R}^{N}} \cal{L}\rho(x)dx \nonumber \\
 & \leq & \cal{L} \| \rho \|_{L^1(R^N)}, 
\end{eqnarray}
where
 $\cal{L}$ is the integral of the continuous function $f^{-1}$  in the (finite) interval
 $ [0, M ]$.

 Finally  let 
 \begin{equation}
 G_3(x) := h(x) f(u(x)) \rho(x) 
 \label{IntegF3} 
 \end{equation}  
denote the integrand of $F_u(u). $ 
 Then
 \[
 |G_3(x) | \leq  M    \| h  \|_{\infty}  \rho(x) 
 \]
and

\begin{eqnarray} \label{boundF3}
|F_{3}(u)|&\leq& \int_{\mathbb{R}^{N}}   M    \| h  \|_{\infty}      \rho(x)dx   \nonumber \\
&\leq&  M  \| h  \|_{\infty}    \| \rho \|_{L^1(R^N)}.
\end{eqnarray}
 \qed



\begin{teo}
  Suppose $f$ satisfies the same hypotheses of Proposition \ref{lower}.
  Then the functional
given in \eqref{Lyap_func}  is continuous in the topology of
$C_{\rho}(\mathbb{R}^{N})$.\label{Contínuo}
\end{teo}
\dem  Write $F(u)=F_{1}(u)+F_{2}(u)-F_{3}(u)$ 
 as in the proof of the Proposition \ref{lower}.

Let $u_n$ be a sequence of functions converging to $u$ in $C_{\rho}(\mathbb{R}^{N})$.

 Let also

 \[  {\displaystyle G_1(x,y), G_2(x), G_3(x)}  \quad    \textrm{ as in} \quad
  \eqref{IntegF1}, \eqref{IntegF2}, \eqref{IntegF3} \quad \textrm{ and } \]
	
	\[  {\displaystyle   G^n_1(x,y), G^n_2(x), G^n_3(x)}  \quad    \textrm{ as in} \quad
  \eqref{IntegF1}, \eqref{IntegF2}, \eqref{IntegF3} \quad \textrm{ with  } u  \textrm{ replaced by } u_n. \]

 Then \[ F_{1} (u_n) =  
 -\frac{1}{2} \int_{\mathbb{R}^{N}} \int_{\mathbb{R}^{N}}\  G^n_1(x,y) dy dx  
 \]
$$
F_{2}(u_n)=\int_{\mathbb{R}^{N}} G^n_2(x) dx 
$$
and
$$
F_{3}(u_n)=\int_{\mathbb{R}^{N}} G^n_3(x) dx.
$$

 By  \eqref{IntegF1}, \eqref{IntegF2}, \eqref{IntegF3}  and  \eqref{boundF1}, \eqref{boundF2}, \eqref{boundF3};
 the integrands $  G^n_1(x,y), G^n_2(x), G^n_3(x) $ are all bounded by integrable functions independent of 
$n$.  Also  from the pointwise convergence of $u_n$ to $u$ and the continuity  of the functions $f, \rho $ and $h$,
it follows that  $  G^n_1(x,y) \to G_1(x,y)  , G^n_2(x) \to G_2(x)$ and $ G^n_3(x) \to G_3(x) $, for all
 $x,y \in \R^N$. 

 Therefore, $F (u_n) \to F (u) $, by Lebesgue Dominated Convergence Theorem

\noindent This completes the proof. \qed

\begin{teo}
Suppose that $f$ satisfies the same hypotheses  of Proposition \ref{lower} and that $ | f'(x) | \leq( |x|+c)\rho^{3}(x)$, for all $x\in \mathbb{R}^{N}$ and some positive constant $c$. Let $u(\cdot,t)$ be a solutions of (\ref{eq_nonloc}). Then $F(u(\cdot,t))$ is differentiable with respect to $t$ and
$$
\frac{dF}{dt}=-\int_{\mathbb{R}^{N}}[-u(x,t)+J{*}_{\rho}(f \circ u)(x,t) +h]^{2}f'(u(x,t))d\mu(x)\leq 0.\label{Lyapunov}
$$
\end{teo}
\dem Let
$$
\varphi(x,s)=-\frac{1}{2}f(u(x,s))\int_{\mathbb{R}^{N}}J(x-y)f(u(y,s)) \rho(y)dy+\int_{0}^{f(u(x,s))}f^{-1}(r)dr
-h f(u(x,s).
$$
Using the hypotheses on  $f$  and  the fact that$ | f'(x) | \leq( |x|+c)\rho^3(x)$, it is easy to see that  $\|\frac{\partial \varphi (\cdot,s)}{\partial
  s}\|_{L^{1}(\mathbb{R}^{N}, d\mu(x))}<\infty$, for all $s\in \mathbb{R}_{+}$. Hence,
derivating
under the integration sign, we obtain
\begin{eqnarray*}
\frac{d}{dt }F(u(\cdot,t)) &=&
\int_{\mathbb{R}^{N}}[-\frac{1}{2}\frac{\partial f(u(x,t))}{\partial t}
\int_{\mathbb{R}^{N}}J(x-y)f(u(y,t))d\mu(y)\\
&-& \frac{1}{2}f(u(x,t))\int_{\mathbb{R}^{N}}J(x-y)\frac{\partial f(u(y,t))}{\partial t}d\mu(y) \\
&+&f^{-1}(f(u(x,t)))\frac{\partial f(u(x,t))}{\partial t}-h\frac{\partial f(u(x,t))}{\partial t}]d\mu(x)\\
&=&-\frac{1}{2}\int_{\mathbb{R}^{N}}\int_{\mathbb{R}^{N}}J(x-y)f(u(y,t))\frac{\partial f(u(x,t))}{\partial t}d\mu(y)d\mu(x)\\
&-& \frac{1}{2}\int_{\mathbb{R}^{N}}\int_{\mathbb{R}^{N}}J(x-y)f(u(x,t))\frac{\partial f(u(y,t))}{\partial t}d\mu(y)d\mu(x)\\
&+& \int_{\mathbb{R}^{N}}[u(x,t)-h]\frac{\partial f(u(x,t))}{\partial
t}d\mu(x).
\end{eqnarray*}
Since
$$
\int_{\mathbb{R}^{N}}\int_{\mathbb{R}^{N}}J(x-y)f(u(y,t))\frac{\partial
f(u(x,t))}{\partial t}d\mu(y)d\mu(x)\\
=\int_{\mathbb{R}^{N}}\int_{\mathbb{R}^{N}}J(x-y)f(u(x,t))\frac{\partial
f(u(y,t))}{\partial t}d\mu(y)d\mu(x)\\,
$$
It follows that

\begin{eqnarray*}
\frac{d}{dt }F(u(\cdot,t)) &=& -\int_{\mathbb{R}^{N}}\int_{\mathbb{R}^{N}}J(x-y)f(u(y,t))\frac{\partial f(u(x,t))}{\partial t}d\mu(y)d\mu(x)\\\\
&+& \int_{\mathbb{R}^{N}}[u(x,t)-h]\frac{\partial f(u(x,t))}{\partial t} d\mu(x)\\
&=&-\int_{\mathbb{R}^{N}}[-u(x,t)+\int_{\mathbb{R}^{N}}J(x-y)f(u(y,t))d\mu(y) +h]\frac{\partial f(u(x,t))}{\partial t}d\mu(x)\\
&=&-\int_{\mathbb{R}^{N}}[-u(x,t)+J{*}_{\rho}(f \circ u)(x,t) +h]\frac{\partial f(u(x,t))}{\partial t}d\mu(x)\\
&=&-\int_{\mathbb{R}^{N}}[-u(x,t)+J{*}_{\rho}(f \circ u)(x,t) +h]f'(u(x,t))\frac{\partial u(x,t)}{\partial t}d\mu(x)\\
&=&-\int_{\mathbb{R}^{N}}[-u(x,t)+J{*}_{\rho}(f \circ u)(x,t)
+h]^{2}f'(u(x,t))d\mu(x).\\
\end{eqnarray*}
Using that $f$ is strictly increasing,  the result follows.\qed


\begin{Remark}
From Theorem \ref{Lyapunov} follows that, if $F(T(t)u_0)=F(u_0)$ for $t\in \mathbb{R}$, then $u_0$ is an equilibrium point for $T(t)$.\label{Remarkequilibrium}
\end{Remark}

\subsection{Gradient property}

We recall that a semigroup, $T(t)$, is {\em gradient} if each
bounded positive orbit is precompact and there exists a continuous Lyapunov
Functional for $T(t)$, (see \cite{Hale}).

\begin{prop} Assume the same  hypotheses from Theorems \ref{Lyapunov} and \ref{Teorema 3.3}. Then the flow generated by
equation (\ref{eq_nonloc}) is gradient.\label{Gradiente}
\end{prop}
\dem The precompacity
of the orbits follows from existence of the global attractor. From Proposition \ref{lower}, Theorem \ref{Contínuo}, Theorem \ref{Lyapunov} and Remark \ref{Remarkequilibrium} follows that the functional given in \eqref{Lyap_func} is a continuous Lyapunov functional.
\qed
\\

As consequence of the Proposition \ref{Gradiente} we have the convergence of the solutions of (\ref{eq_nonloc}) to  the equilibrium point set of T(t) (see \cite{Hale} - Lemma 3.8.2)

\begin{coro}  For any $u\in C_{\rho}(\R )$,  the $\omega$-limit set, $\omega(u)$, of $u$ under $T(t)$ belongs to $E$. Analogously the $\alpha$-limit set, $\alpha(u)$, of $u$ under $T(t)$ belongs to $E$.\label{limit_set}
\end{coro}

Also as a consequence of the Proposition \ref{Gradiente}, we have that the global attractor given in the Theorem \ref{Teorema 3.3} has the following characterization (see \cite{Hale} - Theorem
 3.8.5).
 \begin{teo} Under the same hypotheses from Theorem \ref{Lyapunov},  the attractor
$\cal A$ is  the unstable set of the
equilibrium point set of $T(t)$, that is,
$$
{\cal A} =W^{u}(E),
$$
where $E=\{u\in B_{\rho}(0,R) \,\, : \,\,
u(x)=J{*}_{\rho}(f\circ u)(x)+h\}$.
\end{teo}

\dem Let $u\in \cal A$. Then, there exists a complete orbit through $u$ which is 
contained in  $ \cal A$. Since $A$ is compact,  the   $\alpha$-limit set, $\alpha(u)$, of $u$ under $T(t)$  is nonempty.  By Lemma \ref{limit_set} it belongs to 
$E$ and, therefore,  $u\in W^{u}(E)$.

 Conversely, suppose   $u\in W^{u}(E)$  and  let  $E^{\delta}$ denote  a $\delta$-neighborhood of $E$.  Then, for any $\delta>0$, there exists $\overline{t}$
such that 
  $
T(-t)u\in E^{\delta},
$ 
for any $t \geq \overline{t}$.  Thus,   $u\in  T(t) (E^{\delta})$, for any 
 $t \geq \overline{t}$,  It follows that $u$  is arbitrarily close to ${\cal A}$, so it must belong to $ {\cal A}$.
 
%
%
This concludes the proof. \qed

\section{An example} \label{example}

Motivated by the example given in \cite{Silva3},  we consider the one dimensional case of 
(\ref{eq_nonloc}), with  $f(x)=(1+e^{-x})^{-1}$,
$$
J(x)=\left\{\begin{array}{ccccc}
e^{\frac{-1}{1-x^{2}}},  \,\,if \,\, |x|< 1,\\
                      0,\,\, if \,\,  |x| \geq  1,\\
\end{array}
\right.
$$
that is, we consider the equation
\begin{equation}
\frac{\partial u(x,t)}{\partial t}=-u(x,t)+
\int_{x-1}^{x+1}e^{\frac{-1}{1-(x-y)^{2}}}(1+e^{-u(y)})^{-1} \sqrt{(1+y^{2})^{-1}}dy+ h.
\label{Ex1}
\end{equation}

It  is easy to see that the function $J$ meet all the hypotheses assumed in introduction, that is, $J$ is an even non negative function of class $C^1(\mathbb{R})$. Furthermore, we have:

\begin{Remark}  The function $f$ satisfies the condition  (\ref{dissip1}), with
$\eta=1$ and $ K=\frac{1}{2}$. 

In fact, since $f'(x)=(1+e^{-x})^{-2}e^{-x}>0$, it follows that $1< (1+e^{-x})^{2}\leq 4$,
$\forall \,\, x \in \mathbb{R}^{N}$. Thus
\begin{equation}
\frac{1}{4}\leq (1+e^{-x})^{-2} <1.\label{I}
\end{equation}
Then, since $f''(x)=2(1+e^{-x})^{-3}e^{-2x}-(1+e^{-x})^{-2}e^{-x}$, follows that $|f''(x)| < 3$,
$\forall \,\, x \in \mathbb{R}^{N}$. Hence $f'$ is locally Lipschitz.
Furthermore, follows from (\ref{I}) that
$$
|f(x)-f(y)|=|(1+e^{-x})^{-1} -(1+e^{-y})^{-1}|\leq |x-y|.
$$
In particular, using that $f(0)=\frac{1}{2}$, results
$$
|f(x)|\leq |x|+\frac{1}{2}, \,\, \forall \, x\in \mathbb{R}^{N}.
$$
\end{Remark}

\begin{Remark}
With $\rho(x)=\sqrt{(1+x^{2})^{-1}}$, the hypothesis that $\int_{\mathbb{R}^{N}} \rho(x)dx < \infty$ is
easily verified and $|\rho(x)|\leq 1$, for all $x\in \mathbb{R}$. Furthermore, we also have 
$$
f'(x)=(1+e^{-x})^{-2}e^{-x} \leq (1 +  | x| )  (1+x^{2})\sqrt{(1+x^{2})^{-1}}  = (1 +  | x |) \rho^{3}(x).
$$
\end{Remark}

\begin{Remark}The  hypotheses in the Theorem \ref{Lyapunov} are also satisfied.

In fact, note that $0<|(1+e^{-x})^{-1}|<1$ and $f^{-1}(x)=-\ln(\frac{1-x}{x})$.
Thus it is easy to see that, for $0\leq s\leq 1$,
$$
\left|\int_{0}^{s}-\ln(\frac{1-x}{x})dx\right|\leq \ln 2.
$$
\end{Remark}



Therefore the results of the preview sections are valid for the flow
generated by equation (\ref{Ex1}).





%

\end{document}